\documentclass[letterpaper, 10 pt, conference]{ieeeconf}  

\IEEEoverridecommandlockouts                              
\usepackage{graphics} 
\usepackage{amsmath,amsfonts} 
\usepackage{tikz}
\usepackage{algorithmic,algorithm}
\usepackage{xcolor,colortbl}
\usepackage{siunitx}

\newtheorem{remark}{Remark} 

\newcommand\R{\mathbb{R}}
\newcommand\D{\mathrm{d}}

\title{\LARGE \bf Controlling the bottom topography of a microalgal pond to optimize productivity  }

\author{{\centering Olivier Bernard$^\dagger$
\thanks{$^\dagger$Universit\'e Nice C\^ote d'Azur, Inria BIOCORE, BP93, 06902 Sophia-Antipolis Cedex, France ({\tt\small olivier.bernard@inria.fr})},
Liudi Lu$^\star$,
Jacques Sainte Marie$^\star$,
Julien Salomon$^\star$
\thanks{$^\star$INRIA  Paris,  ANGE  Project-Team,  75589  Paris  Cedex  12,  France   and Sorbonne Universit\'e, CNRS, Laboratoire Jacques-Louis Lions, 75005 Paris, France ({\tt\small liudi.lu@inria.fr}, {\tt\small Jacques.Sainte-Marie@inria.fr}, {\tt\small julien.salomon@inria.fr})
}
}
}

\begin{document}

\maketitle
\thispagestyle{empty}
\pagestyle{empty}

\begin{abstract}
We present a coupled model describing  growth  of microalgae in a raceway cultivation process, accounting for hydrodynamics. Our approach combines a biological model (based on the Han model) and  shallow water dynamics equations that model the fluid into the raceway. We then describe an optimization procedure dealing with the topography to maximize the biomass production over one  cycle (one lap of the raceway).
The results show that non-flat topographies enhance microalgal productivity.
\end{abstract}

\section{INTRODUCTION}

Research on production of biotechnological microalgae has been booming in recent years. The potential of these techniques finds interests for the cosmetics, pharmaceutical, food and - in the longer term - green chemistry and energy applications~\cite{Wijffels2010}. Production is carried out in biophotoreactors that often take the form of a raceway, i.e. a circular basin exposed to solar radiation where the water is set in motion by a paddle wheel~\cite{Chiaramonti2013}. On top of homogenizing the medium for ensuring an equidistribution of the nutrients necessary for algal growth, the main interest of mixing is to guarantee that each cell will have regularly access to light~\cite{Demory2018}. The algae are regularly harvested, and their concentration is maintained around an optimal value~\cite{Munoz-Tamayo2013,Posten2016}. The algal concentration stays generally below 1\%~\cite{Bernard2015}. Above this value, the light extinction is so high that a large fraction of the population is in the dark and does not grow anymore. 

Many phenomena have to be taken into account to represent the entire photoproduction process. In this paper, we develop a coupled model to describe the growth of algae in a raceway, accounting for the light that they receive. More precisely, prolonging the study of~\cite{Bernard2013} with a simpler model, we combine a photosynthesis model, the Han model, and a hydrodynamic model based on the Saint-Venant equations for shallow water flows. This approach enables us to formulate an optimization problem where the raceway topography is designed to maximize the productivity. For this problem, we present an adjoint-based optimization scheme  including constraints to respect the shallow water regime. Our numerical tests show that non-trivial topographies can be obtained. 

The outline of the paper is as follows: in Section 2, we present the  biological and hydrodynamic models underlying our coupled model. In Section 3, we describe the optimization problem and a corresponding  numerical optimization procedure. Section 4 is devoted to the numerical results obtained with our approach. We conclude in Section 5 with some perspectives opened by this work.

\section{COUPLING HYDRODYNAMICAL AND BIOLOGICAL MODELS}
Our approach is based on a coupling between the hydrodynamic behavior of the particles and the evolution of the photosystems driven by the light intensity they received when traveling across the raceway pond.

\subsection{Modelling the photosystems dynamics}
The photosystems we consider are cell units that harvest photons and transfer their energy to metabolism.  The photosystems dynamics can be described by the Han model~\cite{Han2001}. In this compartmental model, the photosystems can be described by three different states: open and ready to harvest a photon ($A$), closed while processing the absorbed photon energy ($B$), or inhibited if several photons have been absorbed simultaneously ($C$).
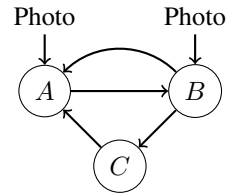
\begin{figure}[thpb]
\begin{center}
\begin{tikzpicture}
\node at (0,0)[circle,draw,scale=1](C){$C$};
\node at (-1,1)[circle,draw,scale=1](A){$A$};
\node at (1,1)[circle,draw,scale=1](B){$B$};
\node at (-1,2) (Pho1){Photo};
\node at (1,2) (Pho2){Photo};
\draw [->,thick](A)--(B);
\draw [->,thick](B)--(C);
\draw [->,thick](C)--(A);
\draw (A) edge[<-,thick,bend left =45] (B);
\draw [<-,thick](A)--(Pho1);
\draw [<-,thick](B)--(Pho2);
\end{tikzpicture}
\end{center}
\caption{Scheme of the Han model, representing the probability to go from one state to another, as a function of the photon flux density.} 
\end{figure}

Their evolution satisfy the following dynamical system
\begin{equation*}
\left\{
\begin{array}{lr}
\dot{A} = -\sigma I A + \frac B{\tau},\\
\dot{B} =  \sigma I A - \frac B{\tau} + k_rC - k_d\sigma I B,\\
\dot{C} = -k_rC + k_d \sigma I B.
\end{array}
\right.
\end{equation*}
Here $A, B$ and $C$ are the relative frequencies of the three possible states with 
\begin{equation}\label{abc}
A+B+C=1,
\end{equation}
and $I$ is the photon flux density, a continuous time-varying signal. Besides, $\sigma$ stands for the specific photon absorption, $\tau$ is the turnover rate, $k_r$ represents the photosystems repair rate and $k_d$ is the damage rate. Following~\cite{Lamare2018} and using~\eqref{abc}, we reduce this system to one single evolution equation:
\begin{equation}\label{evolC}
\dot{C} = -\alpha(I) C + \beta(I),
\end{equation}
where
\begin{equation*}
\begin{split}
\alpha(I) = k_d\tau \frac{(\sigma I)^2}{\tau \sigma I+1} + k_r,\quad  \beta(I) = \alpha(I)-k_r.
\end{split}
\end{equation*}
The net specific growth rate is obtained by balancing photosynthesis and respiration, which gives
\begin{equation}\label{mu}
\mu(C,I) = -\gamma(I)C + \zeta(I),
\end{equation}
where
\begin{equation*}
\begin{split}
\gamma(I)  = \frac{k\sigma I}{\tau \sigma I+1},\quad \zeta(I) = \gamma(I) -R.
\end{split}
\end{equation*}
Here $k$ is a factor that links received energy and growth rate. The term $R$ represents the respiration rate. 

\begin{remark}
The dynamics of the biomass $X$ is derived from $\mu$:
\begin{equation}\label{evolX}
\dot{X} = \mu(C,I)X.
\end{equation}    
\end{remark}

\subsection{Steady 1D Shallow Water Equation}
The shallow water equations are one of the most popular model for describing geophysical flows, which is derived from the free surface incompressible Navier-Stokes equations (see for instance~\cite{Gerbeau2001}).

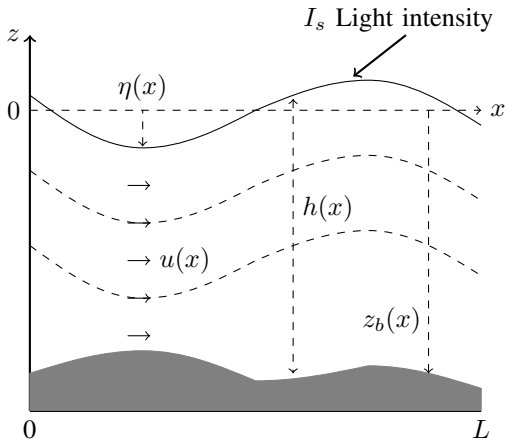
\begin{figure}[thpb]
\begin{center}
\begin{tikzpicture}
\draw [thick](0,-4) -- (6,-4);
\draw (0,-4) node[anchor=north] {0};
\draw (6,-4) node[anchor=north] {$L$};
\draw[thick, ->] (0,-4) -- (0,1) node[anchor=east] {$z$};
\draw (0,0) node[anchor=east] {$0$};
\draw[dashed,->] (0,0) -- (6,0) node[anchor=west]  {$x$};
\draw [dashed, <-](1.5,-0.48) -- (1.5,0);
\draw (1.5,0.3) node{$\eta(x)$};
\draw[thick,<-] (4.3,0.45) -- (5,1);
\draw (4.9,1.2) node {$I_s$ Light intensity};
\draw (0,0.2) sin (1.5,-0.5) cos (3,0) sin (4.5,0.4) cos (6,-0.2);
\draw [dashed] (0,-0.8) sin (1.5,-1.5) cos (3,-1) sin (4.5,-0.6) cos (6,-1.2);
\draw [dashed] (0,-1.8) sin (1.5,-2.5) cos (3,-2) sin (4.5,-1.6) cos (6,-2.2);
\draw[fill, gray] (0,-4) -- (0,-3.5) sin (1.5,-3.2) cos (3,-3.6) cos (4.5,-3.4) cos (6,-3.7) -- (6,-4);
\draw [dashed,->](5.3,0) -- (5.3,-3.5);
\draw (4.3,-2.8) node[anchor=west] {$z_b(x)$};
\draw[dashed,<->] (3.5,0.15) -- (3.5,-3.5);
\draw (3.5,-1.3) node[anchor=west]{$h(x)$};
\draw [->](1.3,-2) -- (1.6,-2) node[anchor=west]{$u(x)$};
\draw [->](1.3,-1.5) -- (1.6,-1.5);
\draw [->](1.3,-2.5) -- (1.6,-2.5);
\draw [->](1.3,-1) -- (1.6,-1);
\draw [->](1.3,-3) -- (1.6,-3);
\end{tikzpicture}
\end{center}
\caption{Representation of the hydrodynamic model.}
\label{physic}
\end{figure}

We assume that the hydrodynamics in the raceway pond corresponds to a laminar flow where the viscosity is neglected. As a consequence, it can be described by the steady Shallow Water equation which is given by
\begin{align}
&\partial_x(hu) = 0, \label{1_eq} \\ 
&\partial_x(hu^2+g\frac{h^2}2) = -g h \partial_x z_b, \label{2_eq}
\end{align}
where $u$ is the horizontal averaged velocity of the water, $h$ is the water elevation, the constant $g$ stands for the gravitational acceleration, and $z_b$ defines the topography of the raceway pond. The free surface $\eta$ is given by $\eta = h + z_b$. This system is presented in Fig~\ref{physic}.
The $z$ axis represents the vertical direction and the $x$ axis represents the horizontal direction. Besides, $I_s$ represents the light intensity at the free surface (assumed to be constant).

Considering Equations~\eqref{1_eq}-\eqref{2_eq}, and  integrating the former, we get
\begin{equation}\label{cst1}
hu=Q_0, 
\end{equation}
for a fixed positive constant $Q_0$. On the other hand, since $h u$ is a constant, Equation~\eqref{2_eq} can be rewritten as
\begin{equation*}
hu\partial_x u + h\partial_x g h + h\partial_x g z_b = 0.
\end{equation*}
Assuming that $h>0$ and using~\eqref{cst1}, we can eliminate $u$ to get
\begin{equation*}
\partial_x ( \frac{Q_0^2}{2h^2} + g (h + z_b) ) = 0.
\end{equation*}
Then for two given constants $h(0), z_b(0)\in \R$, we have for all $x\in [0,L]$
\begin{equation}\label{cst2}
\frac{Q_0^2}{2h^2} + g (h + z_b) = \frac{Q_0^2}{2h^2(0)} + g (h(0) + z_b(0)) =: M_0. 
\end{equation} 
Finally we find that the topography $z_b$ satisfies
\begin{equation}\label{zb}
z_b = \frac{M_0}g - \frac{Q_0^2}{2gh^2} - h.
\end{equation}

\begin{remark}
This model is relevant if the flow is subcritical and the water elevation is positive, meaning that the following relation shall hold
\begin{equation*}
0 < \frac{u(x)}{\sqrt{gh(x)}} < 1 \qquad x\in(0,L),
\end{equation*}
which can also be written as 
\begin{equation}\label{subcrit}
h(x) > (\frac{Q^2_0}{g})^{\frac 13} \qquad x\in(0,L).
\end{equation}
This leads us to introduce the threshold value $h_c := (\frac{Q^2_0}{g})^{\frac 13}$, which guarantees that the system is in a shallow water regime. Note that $h_c$ is the threshold value of $h$ for which the Froude number equals to one. 
\end{remark}

From \eqref{zb}, $h$ is the solution of a third degree equation. Given a smooth topography $z_b$, there exists a unique positive smooth solution of \eqref{zb} which satisfies the subcritical flow condition (see \cite[Lemma 1]{Michel-Dansac2016}).

\subsection{Lagrangian trajectories of the algae and captured light intensity}
Let $z$ be the vertical position of a particle in the raceway. We first determine the Lagrangian trajectory of an algal cell which starts at a given position $z(0)$ at time 0. 
    
From the incompressibility of the flow, we have 
\begin{equation*}
\nabla \cdot \underline{\textbf{u}} = 0,
\end{equation*}
with $\underline{\textbf{u}} = (u(x),w(x,z))$. Here $u(x)$ is the horizontal velocity and $w(x,z)$ is the vertical velocity. This implies that
\begin{equation}\label{incompress}
\partial_x u + \partial_z w = 0.
\end{equation}
Integrating~\eqref{incompress} from $z_b$ to $z$ gives:
\begin{align*}
0&=\int_{z_b}^z \big(\partial_x u(x) + \partial_z w(x,z) \big)\D z, \\
& = \partial_x \int_{z_b}^z  u(x) \D z +\int_{z_b}^z \partial_z w(x,z) \D z,\\
& =  \partial_x \big((z-z_b)u(x)\big) + w(x,z) - w(x,z_b), \\
& = (z - z_b) \partial_x u(x) - u(x)\partial_x z_b + w(x,z).
\end{align*}
where we have used the non-penetration condition $w(x,z_b)=0$. It then follows from~\eqref{zb} that
\begin{equation*}
w(x , z) = (\frac{M_0}g - \frac{3u^2(x)}{2g} - z)u'(x).
\end{equation*}
The Lagrangian trajectory is consequently characterized by the system
\begin{equation}\label{Lagtraj}
\begin{pmatrix}
\dot{x}(t)\\
\dot{z}(t)
\end{pmatrix}
= 
\begin{pmatrix}
u(x(t))\\
w(x(t),z(t)).   
\end{pmatrix}
\end{equation}

To obtain the light intensity observed on this trajectory, we assume the turbidity to be constant over the considered time scale and use Beer-Lambert law:
\begin{equation}\label{Beer}
I(x , z) = I_s\exp{(-\varepsilon (\eta(x) - z))},
\end{equation}
where $\varepsilon$ is the light extinction coefficient. In doing so, we suppose that the system is perfectly mixed so that the concentration of the biomass $X$ defined in~\eqref{evolX} is homogeneous and $\varepsilon$ is constant. 

\begin{remark}
We see that in this approach, the light intensity $I$ couples the hydrodynamic model and Han model: the trajectories of the algae define the received light intensity, which is used in the photosystem dynamics.
\end{remark}
        
\addtolength{\textheight}{-3cm}   

\section{OPTIMIZATION PROBLEM}
In this section, we define the optimization problem associated with the biological-hydrodynamical model. We assume a constant volume for the raceway system. We first introduce our procedure in the case of one single layer, and then extend this to a multiple layers system. For simplicity, we omit $t$ in the notations.

\subsection{One layer problem}

\subsubsection{Optimization functional}
The average net specific growth rate for one trajectory is defined by
\begin{equation}
\bar \mu_1:= \frac 1T \int_0^T \mu( C,I(x,z) ) \D t.
\end{equation}

As we have mentioned in the previous section, in the subcritical case, a given topography $z_b$ corresponds to a unique water height $h$. On the other hand, the volume of our system is defined by
\begin{equation}\label{vol}
V = \int_0^L h(x) \D x.
\end{equation}
Therefore, we choose to parameterize $h$ by a vector $a\in \R^{N}$, which will be the variable to be optimized, in order to handle the volume of our system. 
Given a vector $a$ and the associated $h$, the optimal topography can be obtained by using~\eqref{zb}.

Our goal here is to optimize the topography to maximize $\bar \mu_1$. In this way, we define the functional with the help of~\eqref{mu}
\begin{equation*}
\bar \mu_1(a) = \frac 1T\int_0^T -\gamma(I( x,z;a ))C + \zeta(I( x,z;a )) \D t,
\end{equation*}
where $C, x, z$ satisfy
\begin{equation}\label{cont}
\left\{
\begin{array}{lr}
\dot{C} = -\alpha(I(x,z;a)) C + \beta(I(x,z;a))\\
\dot{x} = u(x ; a)\\
\dot{z} = w(x,z;a).
\end{array}
\right.
\end{equation}
The optimal control problem then reads:

\textit{Find $a^\star$ solving the maximization problem:}
\begin{equation}\label{Jopt}
\max_{a\in \R^{N}} \bar \mu_1(a).
\end{equation}

\subsubsection{Optimality System}
Define the Lagrangian of Problem~\eqref{Jopt} by
\begin{align*}
\mathcal{L}(C,z,&x,p_1,p_2,p_3,a) \\
&:=\frac 1T \int_0^T \big( -\gamma(I(x,z;a))C + \zeta(I(x,z;a)) \big) \D t\\
&- \int_0^T p_1 \big( \dot{C} + \alpha( I(x,z;a) ) C - \beta( I(x,z;a) ) \big) \D t\\
&- \int_0^T p_2 \big( \dot{z} - w(x,z;a) \big) \D t- \int_0^T p_3 \big( \dot{x} - u(x;a) \big) \D t.
\end{align*}
where $p_1, p_2$ and $p_3$ are the Lagrange multipliers associated with the constraints~\eqref{cont}.

The optimality system is obtained by cancelling all the partial derivatives of $\mathcal{L}$. Differentiating $\mathcal{L}$ with respect to $p_1, p_2, p_3$ and equating the resulting terms to zero gives the corrected model equations~\eqref{cont}. Integrating the terms $\int p_1\dot{C}\D t$, $\int p_2\dot{z}\D t$ and $\int p_3\dot{x}\D t$ on the interval $[0,T]$ by parts and differentiating $\mathcal{L}$ with respect to $C,z,x,C(T),z(T),x(T)$ gives rise to
\begin{equation}\label{adjoint}
\left\{
\begin{array}{lr}
\partial_C\mathcal{L} = -\frac 1T \gamma(I(x,z;a)) + \dot{p}_1 -\alpha(I(x,z;a))p_1\\
\begin{aligned}
\partial_z\mathcal{L} &= \Big(\frac 1T\big( -\gamma'( I(x,z;a) )C + \zeta'( I(x,z;a) ) \big) \\
&+p_1\big(-\alpha'(I(x,z;a))C + \beta'(I(x,z;a))\big)\Big) \\
&\partial_zI(x,z;a) + \dot{p}_2 +p_2\partial_z w(x,z;a)
\end{aligned}\\
\begin{aligned}
\partial_x\mathcal{L} &= \Big(\frac 1T\big( -\gamma'( I(x,z;a) )C + \zeta'( I(x,z;a) ) \big) \\
&+p_1\big(-\alpha'(I(x,z;a))C + \beta'(I(x,z;a))\big)\Big) \\
&\partial_xI(x,z;a)+ p_2\partial_x w(x,z;a) + \dot{p}_3 +p_3\partial_x u(x;a)
\end{aligned}\\
\partial_{C(T)}\mathcal{L}  = p_1(T)\\
\partial_{z(T)}\mathcal{L}  = p_2(T)\\
\partial_{x(T)}\mathcal{L}  = p_3(T).
\end{array}
\right.
\end{equation}
Given a vector $a$, let us still denote by $C, x, z, p_1, p_2, p_3$ the corresponding solutions of~\eqref{cont}  and~\eqref{adjoint}. The gradient $\nabla \bar \mu_1(a)$ is obtained by
\begin{equation*}
\nabla \bar \mu_1(a) = \partial_{a}\mathcal{L},
\end{equation*}
where
\begin{equation}\label{grad}
\begin{split}
\partial_a\mathcal{L} = \frac 1T&\int_0^T \big(-\gamma'(I(x,z;a))C + \zeta'(I(x,z;a)) \big)\\
&\partial_aI(x,z;a)\D t + \int_0^Tp_1\big(-\alpha'(I(x,z;a))C \\
+ &\beta'(I(x,z;a))\big)\partial_aI(x,z;a) \D t\\
+ &\int_0^T p_2\partial_a w(x,z;a) \D t + \int_0^T p_3\partial_a u(x;a) \D t.
\end{split}
\end{equation}

\subsection{Multiple Layers Problem}
We now extend the previous procedure to deal with multiple layers. Let us denote $N_z$ the number of layers and $C_i$ (resp. $z_i$) the photo-inhibition state (resp. the trajectory position) associated with the $i$-th layer. As the optimization functional, we consider the semi-discrete average net specific growth rate over the domain, namely:
\begin{equation}\label{mulJ}
\bar \mu_{N_z}(a) = \frac 1{N_z}\sum_{i=1}^{N_z} \frac 1T \int_0^T\mu(C_i,I( x,z_i; a) ) \D t,
\end{equation}
where $C_i, x, z_i$ verify the constraints~\eqref{cont} for $i=1,\cdots,N_z$.
\begin{remark}
Note that the average net specific growth rate over the domain is defined by:
\begin{equation*}
\bar \mu_{\infty} := \frac 1L\int_0^L\frac 1{h(x)}\int_{z_b(x)}^{\eta(x)} \mu(C, I(x,z)) \D z \D x.
\end{equation*}
Our approach consequently consists in considering a vertical discretization of $\bar \mu_{\infty}$, which gives~\eqref{mulJ}. This discretization should not give rise to any problem and is left to a future contribution.
\end{remark}

Similar computations as in the previous section give rise to $N_z$ systems similar to~\eqref{adjoint}, where $C$ and $z$ are replaced by $C_i$ and $z_i$ respectively. Denoting by $p_{1,i}, p_{2,i}$ the associated Lagrange multipliers, the partial derivatives $\partial_x\mathcal{L}$ is a little different from the previous section, more precisely, this can be computed by: 
\begin{equation*}
\begin{split}
\partial_x\mathcal{L} = &\frac 1{N_z}\sum_{i=1}^{N_z}\frac 1T\big( -\gamma'( I(x,z_i;a) )C_i+ \zeta'( I(x,z_i;a) ) \big)\\
&\partial_xI(x,z_i;a)+\sum_{i=1}^{N_z}p_{1,i}\big(-\alpha'(I(x,z_i;a))C_i \\
&+ \beta'(I(x,z_i;a))\big)\partial_xI(x,z_i;a)\\
&+ \sum_{i=1}^{N_z}p_{2,i}\partial_x w(x,z_i;a)+ \dot{p}_3 +p_3\partial_x u(x;a).
\end{split}
\end{equation*}
Finally, the gradient $\nabla \bar \mu_{N_z}(a)$ is given by 
\begin{equation*}
\begin{split}
\partial_a\mathcal{L} &= \frac 1{N_z}\sum_{i=1}^{N_z}\frac 1T\int_0^T\big( -\gamma'( I(x,z_i;a) )C_i+ \zeta'( I(x,z_i;a) ) \big)\\
&\partial_aI(x,z_i;a)\D t+\sum_{i=1}^{N_z}\int_0^T p_{1,i}\big(-\alpha'(I(x,z_i;a))C_i\\
&+\beta'(I(x,z_i;a))\big)\partial_aI(x,z_i;a)\D t \\
&+\sum_{i=1}^{N_z} \int_0^T p_{2,i}\partial_a w(x,z_i;a)\D t + \int_0^Tp_3\partial_a u(x;a)\D t.
\end{split}
\end{equation*}

\section{NUMERICAL EXPERIMENTS}

\subsection{Numerical method}

\subsubsection{Gradient Algorithm \label{pargrad}}
In order to tackle the optimization problem~\eqref{Jopt}, we consider a gradient-based optimization method. The complete procedure is detailed in Algorithm~\ref{algo}. 
\begin{algorithm}
\begin{algorithmic}
\STATE \textbf{Input}: Tol$>0$, $\rho>0$.
\STATE \textbf{Initial guess}: $a$.
\STATE \textbf{Output}: $a$
\STATE Set $err :=$ Tol$+1$ and define $h$ by \eqref{Fourier} using the input data.
\WHILE{$err>$Tol and $\|h\|_{\infty} >h_c$}
\STATE Compute $u$ by \eqref{cst1}.
\STATE Set $x, z$ as the solutions of the last two equations of \eqref{cont}.
\STATE Compute $I$ by \eqref{Beer}.
\STATE Set $C$ as the solution of the first equation of \eqref{cont}.
\STATE Set $p_1, p_2, p_3$ as the solutions of \eqref{adjoint}.
\STATE Compute the gradient $\nabla \bar \mu_{1}$ by \eqref{grad}.
\STATE $a = a + \rho \nabla \bar \mu_{1}$,
\STATE Set $err := \| \nabla \bar \mu_{1}\|$. 
\ENDWHILE
\end{algorithmic}
\caption{Gradient-based optimization algorithm}
\label{algo}
\end{algorithm}
Note that in addition to a numerical tolerance criterion on the magnitude of the gradient, we have added a constraint on the water height $h$. The latter guarantees that we remain in the framework of subcritical flows~\eqref{subcrit} (and in the range of industrial constraints, see~\cite{Chiaramonti2013}).

A similar algorithm can be considered to tackle the multiple layers Problem~\eqref{mulJ}. Remark that since there is no interaction between layers, the gradient computation can be partially parallelized when computing $z_i, C_i, p_{1,i}, p_{2,i}$.

\subsubsection{Numerical Solvers}
We introduce a supplementary space discretization with respect to $x$ to solve our optimization problem numerically, in this way, we set a time step $\Delta t$, and use Heun scheme to compute the trajectories described by~\eqref{Lagtraj}. The integration is stopped when $x$ reaches $L$. Denoting by $N_T$ the final number of the time steps, we then have $N_T\Delta t\approx T$. 

We use the Heun scheme again for computing $C$ via~\eqref{evolC}. We use a first-discretize-then-optimize strategy, meaning that the Lagrange multipliers $p_1$ (resp. $p_2, p_3$) are also computed by a (backward) Heun's type scheme. Note that this scheme is still explicit since it solves a backward dynamics starting from $p_1(T) = 0$ (resp. $p_2(T) = 0, p_3(T)=0$).

\subsection{Parameter settings}

\subsubsection{Parameterization}
In order to describe the bottom of an optimized raceway pond, we choose to parameterize $h$ by a truncated Fourier series for our numerical tests so that $h(0)=h(L)$, which means we have accomplished one lap. More precisely, $h$ reads:
\begin{equation}\label{Fourier}
h(x) = a_0 + \sum_{n=1}^N a_n \sin(2n\pi \frac xL).
\end{equation}
The parameter to be optimized is the Fourier coefficients $a:=[a_1,\cdots, a_N]$. Note that we choose to fix $a_0$, since it is related to the volume $V$ of the raceway. Indeed, under this parameterization we have
\begin{equation*}
V = \int_0^L h(x) \D x =a_0L.
\end{equation*}
From~\eqref{cst1} and~\eqref{zb}, the velocity $u$ and the topography $z_b$ read also as functions of $a$. Once we find the vector $a$ maximizing the functional $\bar \mu_{N_z}$, we then find the optimal topography of our system.

\subsubsection{Parameters of the system}
The time step is set to $\Delta t=\SI{0.1}{s}$ which corresponds to the numerical convergence of the Heun's scheme, and we take $Q_0 = \SI{0.04}{m^2.s^{-1}}$, $a_0=\SI{0.4}{m}$, $z_b(0)=-\SI{0.4}{m}$ to stay in standard ranges for a raceway. We take here $L=\SI{10}{m}$ and the free-fall acceleration $g=\SI{9.81}{m.s^{-2}}$. The initial state of $C$ is set to be its steady state. All the numerical parameters values for Han's model are taken from~\cite{Grenier2020} and recalled in Tab.~\ref{Tab1}.
\begin{table}[htbp]
\caption{Parameter values for Han Model }
\label{Tab1}
\begin{center}
\begin{tabular}{|c|c|c|}
\hline
$k_r$  & $6.8$ $10^{-3}$ & s$^{-1}$\\
\hline
$k_d$ & $2.99$ $10^{-4}$  & -\\
\hline
$\tau$ & 0.25 & s\\
\hline
$\sigma$ & 0.047 & m$^2$.($\mu$mol)$^{-1}$\\
\hline
$k$  & $8.7$ $10^{-6}$ & -\\
\hline
$R$ &  $1.389$ $10^{-7}$ & s$^{-1}$\\
\hline
\end{tabular}
\end{center}
\end{table}

In order to determine the light extinction coefficient $\varepsilon$, let us assume that only $10\%$ light can be captured by the cells at the bottom of the raceway, meaning that $I_b=0.1I_s$, we choose $I_s=\SI{2050}{\mu mol.m^{-2}.s^{-1}}$ which corresponds to the average light intensity at summer (see~\cite{Grenier2020}). Then $\varepsilon$ can be computed by
\begin{equation*}
\varepsilon = (1/a_0)\ln(I_s/I_b).
\end{equation*}

\subsection{Numerical results}

\subsubsection{Convergence test}
\begin{figure}[htpb]
\centering
\includegraphics[scale=0.25]{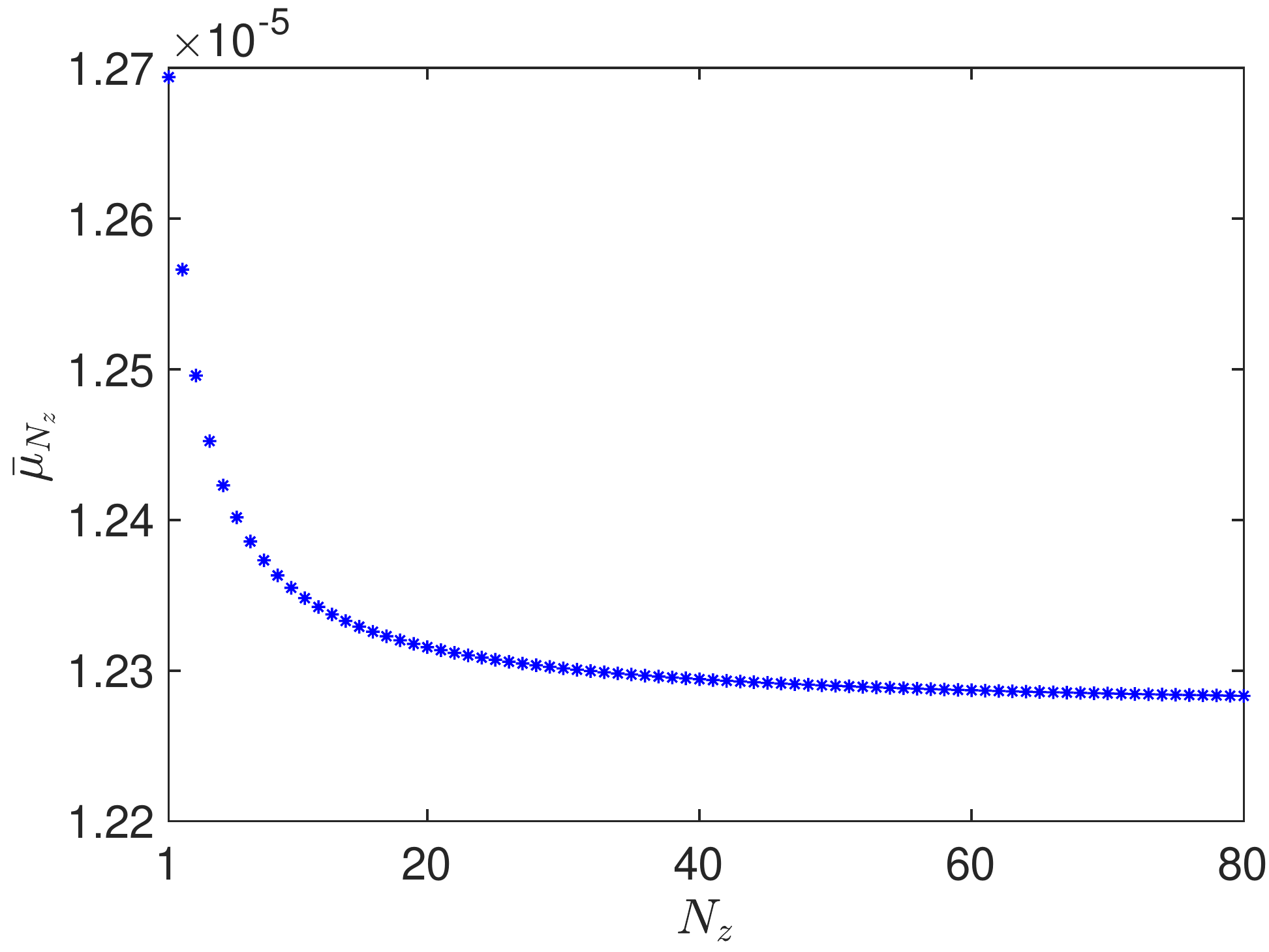}
\caption{The value of $\bar \mu_{N_z}$ for $N_z=[1,80]$.}
\label{figJNz}
\end{figure}

The first test consists in studying the influence of the vertical discretization number $N_z$. We fix $N=5$ and take 100 random initial guesses of $a$. Note that the choice of $a$ should respect the subcritical condition \eqref{subcrit}. For $N_z$ varying from $1$ to $80$, we compute the average value of $\bar \mu_{N_z}$ for each $N_z$. The results are shown in Fig.~\ref{figJNz}.
We observe numerical convergence when $N_z$ grows, showing the convergence towards continuous model. In view of these results, we take $N_z=40$ for the successive studies.


\subsubsection{Optimization example}

We focus now on the shape of the optimal topography. Let us set the numerical tolerance Tol$=10^{-10}$, we choose $N=5$ terms in the truncated Fourier series as an example to show the shape of the optimal topography. As an initial guess, we consider the flat topography, meaning that $a$ is set to 0.
\begin{figure}[htpb]
\centering
\includegraphics[scale=0.3]{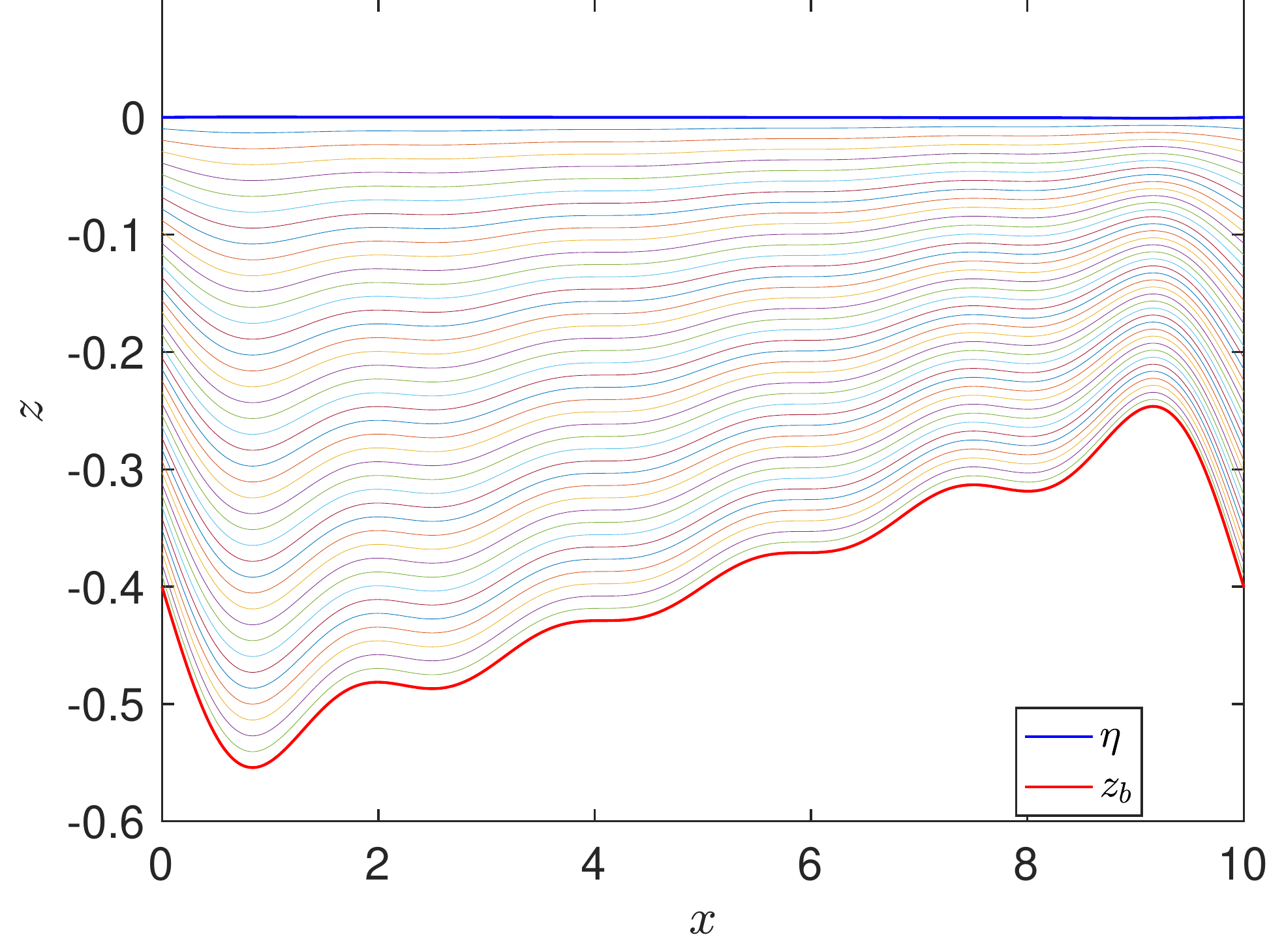}
\includegraphics[scale=0.3]{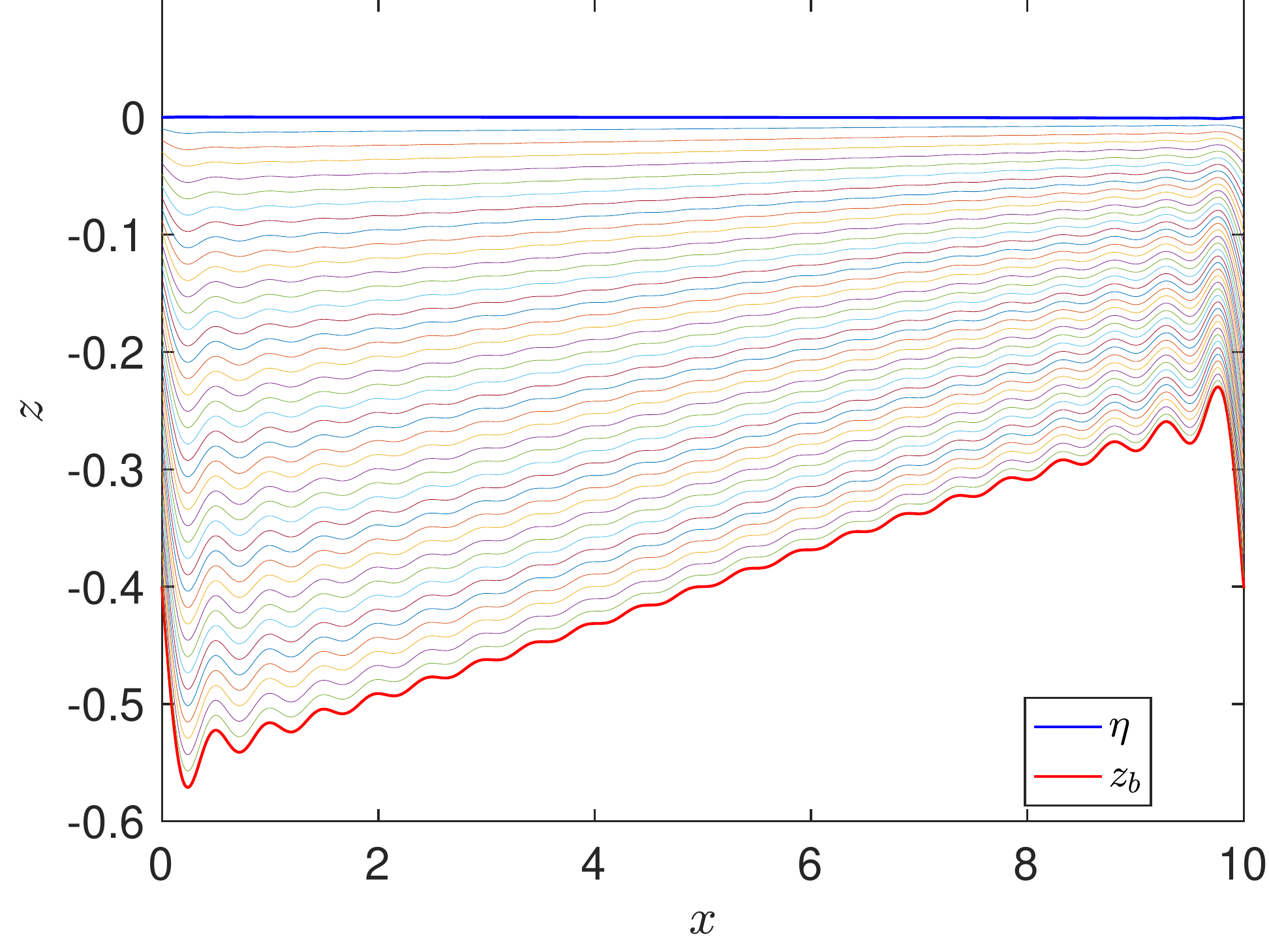}
\caption{The optimal topography at the final iteration for $N=5$ (up) and for $N=20$ (down). The red thick curve represents the topography ($z_b$), the blue thick curve represents the free surface ($\eta$), and all the other curves between represent the trajectories for different layers.}
\label{figN}
\end{figure}
The optimal shape of the topography is shown in Fig.~\ref{figN}, and the $a^\star$ for the final iteration reads $a^\star=[0.1043,0.0503,0.0333,0.0250,0.0201]$. The resulting optimal topography is not flat, which is not a common knowledge in the community.

\subsubsection{Study of the influence of $N$}
The last test is given to study the influence of the order of the truncation of the Fourier series. Set $N=[0,5,10,15,20]$ and keep all the other parameter settings. Tab.~\ref{tabN} shows the optimal value of our functional $\bar \mu_{N_z}$ for different values of $N$. Note that, even the value of $\bar \mu_{N_z}$ is small, it is not the numerical error. There is a slight increase of the optimal value of the functional $\bar \mu_{N_z}$ when $N$ becomes larger. However, corresponding values of $\bar \mu_{N_z}$ remain close to the one associated with a flat topography. As for the optimal shape, for instance, we give in Fig.~\ref{figN} another optimal topography when $N=20$. The optimal topography seems to converge to a slope.
Here, the topography varies significantly
 so that viscosity of the fluid can give rise to turbulence (and irreversibility) which is not taken into account in~(\ref{1_eq}--\ref{2_eq}).

\begin{table}[htpb]
\caption{The value of $\bar \mu_{N_z}$ for different $N$}
\label{tabN}
\begin{center}
\begin{tabular}{|c|c|c|c|c|}
\hline
$N$ & Iter & $\bar \mu_{N_z}$ & $\log_{10}(\|\nabla \bar \mu_{N_z}\|)$\\
\hline
0 & 0 & 1.232270 10$^{-5}$& $-$ \\
\hline
5 & 39 & 1.250805 10$^{-5}$ & -10.077503\\
\hline
10 & 39 & 1.251945 10$^{-5}$ & -10.091920\\
\hline
15 & 39 & 1.252354 10$^{-5}$ & -10.097072\\
\hline
20 & 39 & 1.252565 10$^{-5}$ & -10.099389\\
\hline
\end{tabular}
\end{center}
\end{table}

\section{CONCLUSIONS AND FUTURE WORKS}

A non flat topography slightly enhances the average growth rate. However the gain stays very limited and it is not clear if the practical difficulty to generate such pattern would be compensated by the increase in the process productivity.

Here the class of functions representing the topography has been chosen so that the topographic functions are zero on average. As a consequence, their choice does not affect the liquid volume. Further studies could consist in considering more general class of functions, also considering a non constant volume. Moreover, the optimal algal biomass concentration is itself related to the growth rate (\cite{Martinez201811}) so that the light extinction coefficient is no more a constant, and depends on biomass supported by the parameterization of the system. The computation of the gradient should then be revisited to manage this more complicated case. 


\section{ACKNOWLEDGEMENTS}

This research benefited from the support of the FMJH Program PGMO funded by EDF-THALES-ORANGE.


\bibliographystyle{plain}
\bibliography{auto}  

\end{document}